\input amstex.tex
\input amsppt.sty

\magnification=1100
\topmatter
\title Equivariant spectral flow and a Lefschetz theorem on odd 
dimensional spin manifolds\endtitle
\author Hao Fang\endauthor
\affil Mathematical Sciences Research Institute \endaffil
\address 1000 Centennial Drive, Berkeley, CA 94720 \endaddress
\email haofang{\@}msri.org \endemail
\date May, 2001 \enddate

\rightheadtext\nofrills{EQUIV. SPECTRAL FLOW AND A LEFSCHETZ FORMULA}  

\endtopmatter
\NoBlackBoxes
 
\document 
  
\def\R{\Bbb R}
\def\C{\Bbb C}  

\def\N{\Bbb N}
\def\tensor{\otimes}

\define\coD #1{\nabla_{#1}}
\define\pcoD #1#2{\nabla^{#1}_{#2}}
\define\LC #1{\nabla^{#1}}
\define\dirac {\Delta}
\define\ex {\wedge}
\define\con {\roman{\iota}}
\define\norm #1{\left\| {#1}\right\|}
\define\etabar {\bar \eta  }
\define\half {{1 \over 2}}
\def\e{\epsilon}

\def\lap{\Delta}
\def\im{\sqrt{-1}}
\def\la{\lambda}
\def\th{\theta}
\def\Ga{\Gamma}
\def\FJ{\Cal J}
\def\ind{{\roman {Ind}}}
\def\tr{{\roman {Tr}}}
\def\coker{\roman{coker}}
\def\pf{{\roman {Pf}}}

\def\ch{{\roman {ch}}}
\def\sf{{\roman{sf}}}
\def\spec{{\roman{Spec}}}  

\def\CH{\Cal H}
\def\C{\Bbb C}
\def\R{\Bbb R}

\def\tu{\tilde u}
\def\F{\Cal {F}}
\TagsOnRight

\heading Section 0: Introduction\endheading
  
As one of the most important theories in mathematics, Atiyah-Singer index
theorems have various profound applications and consequences. At the same 
time, there are
several ways to prove these theorems.  Of  particular
interest  is the heat kernel proof, which allows one to obtain
refinements of the index theorems,
i.e.,  the local index theorems for Dirac
operators.  Readers are referred to~\cite{BGV} for a comprehensive treatment 
of the heat
kernel method on even dimensional manifolds.
It is worthwhile to point out here that the heat kernel method also 
lead to direct analytic proofs of the equivariant index theorem for 
Dirac operators on even dimensional spin manifolds. Among the existing 
proofs we list Bismut~\cite{BV}, Berline-Vergne~\cite{BV} and
Lafferty-Yu-Zhang~\cite{LYZ}.

The purpose of this paper is to present a heat kernel proof of 
an equivariant index theorem on odd dimensional spin manifolds, which
is stated for Toeplitz operators.

Recall that Baum-Douglas~\cite{BD} first stated and proved an odd index 
theorem for Toeplitz operators  using the general Atiyah-Singer index theorem
for elliptic pseudo-differential operators. It is known to experts that one 
can give a heat kernel proof of the above mentioned odd index 
theorem. However, let us still  give a brief
description of the basic ideas. The first step is to apply a result
of Booss-Wojciechowski~\cite{BW} to identify the index of the Toeplitz operator
to the spectral flow of a certain family of self-dual elliptic operator 
with positive order. The second step is then to use the well-known
relationship between spectral flows and variations of $\eta$-invariants
to evaluate this spectral flow (cf.~\cite{G}).

Our proof of the equivariant odd index theorem follows the same strategy.
For this purpose we need to introduce a concept of  equivariant  spectral 
flow and establish an equivariant version of the Booss-Wojciechowski theorem
mentioned above. We then extend the relationship between the spectral
flow and variations of $\eta$ invariants to the equivariant 
setting. Finally, we use the local index techniques to evaluate these 
variations.

Among the methods of Bismut~\cite{B}, Berline-Vergne~\cite{BV} and
Lafferty-Yu-Zhang \cite{LYZ}, for simplicity we will follow those
of~\cite{LYZ}
in this paper. There is, however, no  difficulty applying other methods.

Also notice that Dai and Zhang~\cite{DZ} introduced the concept of higher 
spectral flow and gave a heat kernel treatment to the family index 
problem for Toeplitz operators.

This paper is organized as follows. In Section One, we review the basic
definition of the Toeplitz operators associated to Dirac operators
on odd dimensional spin manifolds and prove the equivariant odd
index theorem by using the Baum-Douglas~\cite{BD} trick
and also the general Atiyah-Singer Lefschetz fixed point theorem~\cite{AS}
for elliptic pseudo-differential operators. In Section Two, we introduce the
equivariant spectral flow and prove an equivariant extension of the
Booss-Wojciechowski theorem~\cite{BW}. In Section Three, we establish a
relation between the equivariant spectral flow and the variations of
equivariant $\eta$ invariants. This in turn gives a heat kernel
formula for the equivariant index of the Toeplitz operators. In Section
Four, we evaluate these variations by adopting the local index theorem 
techniques.

The author would like to thank W. Zhang for helpful discussions. He
is also  grateful to X. Dai and B. Cui for their  critical reading
on a preliminary version of this paper. 

\heading Section 1: Toeplitz operators and a Lefschetz fix point theorem
\endheading
 
We begin by fixing notations on  odd dimensional 
Clifford algebras that are used in this paper. From now on, we fix 
$n=2m+1$, where $m$ is a positive integer.

Let $V$ be a $n$-dimensional real vector space associated with a positive
 inner product. Set $e_1, \dots , e_n$ to be an orthonormal basis 
for $V$.
Denote $$T(V)=\R \oplus V \oplus V\otimes V \oplus \cdots,\tag 1.1$$
and $I$ to be the two-sided ideal of $T(V)$ generated by $\{x\otimes 
x+(x,x)1;x\in V\}$.
The Clifford algebra associate to $V$ is defined as $$C(V)=T(V)/I.\tag 1.2$$ 
which is also called as $C(n)$ sometimes. It is clear that 
$\{c_i=e_i I\in C(V)\}$ is the set of generators of $C(V)$ satisfying the 
following relations:
$$c_i c_j+c_j c_i=-2\delta_{ij}.\tag 1.3$$

Define the  chirality operator of $C(n)\otimes \C$ to be
$$\Gamma=(\im)^{m+1} c_1\dots c_n,\tag 1.4$$
which can be checked to be in the center of $C(n)\otimes \C$. It 
is  known  that there is a unique irreducible complex $C(n)$ 
representation $S$ of  dimension $2^m$, such that $\Gamma=Id_{S}$ on $S$.

For future use, we define the symbol map $\sigma: C(n)\to End(\wedge \C
^n)$:
$$\sigma(c_i)= e_i\ex -\con(e_i^*),\tag 1.5$$
which is a complex representation of 
$C(n)$. For $x\in C(n)$ such that  $x$ is not a scalar, we 
have $$\tr_S(x)=-\im (-2\im)^m(\sigma(x)1)_{[n]},\tag 1.6$$
where $(.)_{[d]}$, for any integer $d$, denotes the $d-$dimensional 
part of an exterior form.

We proceed to define the Toeplitz operator on a closed spin manifold.

Through out this paper, we assume $M$ to be a closed (compact, without
boundary), oriented, spin manifold with dimension $n=2m+1$ and a fixed 
spin structure. We also fix a Riemannian metric $g_{TM}$ on $M$.

Let $S(M)$ be the canonical complex spinor 
bundle of $M$,  which  is also a $C(T^*M)$-module.  Let $\LC{TM}$ 
be the canonical Levi-Civita 
connection, which induces a natural connection $\LC{S}$ on $S$. Choose a 
local orthonormal basis 
$e_1,\dots ,e_n$ is  for $TM$, with dual basis $e^1,
\dots ,e^n\in T^*M$. The canonical Dirac operator on $S$ can be defined to be
$$D^S=\sum{c(e^i) {\pcoD {S} {e_i}}}.\tag 1.7$$

It is well known that $D^S$ is a
self-adjoint first-order elliptic differential operator acting on $S(M)$. 
Therefore, there is a spectral decomposition of $\Gamma_{L^2}(S)$ according 
to $D^S$. Denote $L^2_+(S)$ to be
the direct sum of eigenspaces of $D$ associated to nonnegative
eigenvalues, and $P_+$ to be the orthogonal projection operator from
$L^2(S)$ to $L^2_+(S)$. Set $\ P_-=Id-P_+$ and $$P=P_+ -P_-.\tag 1.8$$

Given $\C^N$ a trivial complex vector bundle over $M$ carrying the 
trivial metric and connection,  $D$ and $P$  extend trivially as  
operators acting on $\Gamma(S\tensor\C^N)$.
Let $g:M\to U(N)$ be a smooth map. Then, $g$ extends to an action on 
$S(M)\tensor\C^N$ as $Id_{S(M)}\tensor g$, which is still denoted as $g$ for 
simplicity.
 
\proclaim{Definition 1.1}
Define the Toeplitz operator associated to $D$ and $g$ to be
$$ T_g=(P_+\tensor Id_{\C^N})g(P_+\tensor Id_{\C^N}):
L^2_+(S(M)\tensor\C^N)\to L^2_+(S(M) \tensor\C^N). \tag 1.9$$
\endproclaim

It is a classical fact that $T_g$ is a bounded  
 Fredholm operator between the given Hilbert spaces. Furthermore, if we 
define $\Gamma_\la$ to be the eigenspace of $D$ with eigenvalue $\la$, 
$\Gamma_\la$ is of finite dimension for each $\la$.
 
We then describe the equivariant index problem for Toeplitz operators.

Consider $H$ a compact group of isometries of $M$ preserving the
orientation and spin structure, hence, it also acts on 
$\Gamma(S(M)\tensor\C^N)$. Since the action of $H$ commutes
with the Dirac operator $D$, it  also commutes with $P_{+}$, and $P$. 
Furthermore each $\Gamma_\la$ is $H-$invariant. But to ensure the 
$H$-invariance of Toeplitz operator, we need to make the following assumption on $H$:

\proclaim{Assumption 1.2} For $h\in H$ and any $x\in M$, 
$$g(hx)=g(x).\tag 1.10$$
\endproclaim

As a consequence, 
$$T_g h_{\Gamma(S(M)\tensor\C^N)}=h_{\Gamma(S(M)\tensor\C^N)} T_g.\tag 1.11$$

\proclaim{Definition 1.3}Given $T_g$ and $H$ as above and satisfying (1.10),
the equivariant index of $T_g$, associated with $H$, is defined as the 
following virtual representation of $H$ in $R(H)$, the representation ring of $H$:

$$\ind_H(T_g)=\ker\ T_g-\coker\ T_g.\tag 1.12$$
We also denote, for any $h\in H$,
$$\ind(h, T_g)=\tr(h,\ind_H(T_g)).\tag 1.13$$
\endproclaim
 
Now an application of the general Atiyah-Singer index theorem~\cite{AS} as in 
Baum-Douglas~\cite{BD} gives the following 

\proclaim{Theorem 1.4} For $T_g$ defined as above, let $F_i$'s be 
the fixed, connected sub-manifolds of $M$ under 
the action of any $h\in H$, and  $\nu_i$ be the normal bundle of $F_i$ in 
$TM$, then we have 
$$ \ind(h,T_g)=\sum_i {
({-\im\over{2\pi}})^{m+1-{\dim F_i\over 2}} (\hat A(F) \ch(g)
[\pf (2\sin (\im{(R^\nu(F_i) +\Theta_i)\over 2}))]^{-1}[F_i]}, \tag 1.14$$
where under any local coordinate system,  $\Theta_i$ is the logarithm of the 
Jacobian matrix of $h|\nu_i$, $R^{\nu_i}$
is the curvature matrix of the bundle $\nu_i$, and
$$\ch(g)=\int_0^1 {\tr[g^{-1}dg \exp(u(1-u)(g^{-1} dg)^2)]du}\tag 1.13$$
is the so-called odd Chern character for the differentiable map $g: M\to 
U(N)$. \endproclaim

In Section Four, we will prove a local version of Theorem 1.4.

\heading Section 2: Equivariant spectral flow and equivariant index problem
\endheading

In this section we will introduce the equivariant spectral flow and 
discuss its relation with the equivariant index problem that we have 
set up in the previous section.

Denote $I=[0,1]$. Let $\{D_u\}_{u\in I}$ is a continuous family of self-dual 
elliptic operators of the positive orders on 
the Hilbert space $\CH=L^2(S\otimes \C^N)$. For any fixed $u\in I$, $\spec 
D_u$ is discrete, so we can denote 
the corresponding eigenspace as $\Gamma_{u,\la} $ for any $\la\in \spec D_u$. 
Furthermore, for any 
open $U\in \R$, define $\Gamma_{u,U}=\oplus_{\la\in U} \Gamma_{u,\la}$.

Recall first the usual (scalar) spectral flow according to 
[APS2].

\proclaim{Definition 2.1}
For  $D_u (u\in I)$ a continuous family of self-dual elliptic 
operator of positive order, consider the graph of $\ \spec D_u$:
$${\frak S}=\cup \spec D_u,\tag 2.1$$
which is a closed set of $\R\times I$. Define the spectral flow of 
$\{D_u\}$ to be the intersection number of ${\frak S}$
with the line $\{-\delta\}\times I$ for a sufficiently small positive
$\delta$, which is denoted as $\sf(\{D_u\})$.
\endproclaim

Notice that if both $D_0$ and $D_1$ are invertible, we can simply replace 
$\delta$ in the above definition by $0$.

We then would like to extend this notion to the equivariant case.

Let $H$ be  as in Section One and  $R(H)$ be its representation 
ring. We  further assume that each $D_u$ in the above discussion is 
compatible with the action of $H$. Thus, every $\Ga_{u,\la}$ can be 
viewed as an element of $R(H)$.

We establish the following lemma, which is an extension of 
continuity of the spectrums of the family of self-dual operators.

\proclaim {Lemma 2.2}
Let $\{D_u\}$ be described as above. For a fixed
$u_0\in I$ and any $\la \in \spec D_{u_0}$ with $\dim \Ga_{u_0,\la}=k$, we
can find a positive $\e$ such that for any $u\in I$, $|u-u_0|<\e$, there 
is an open set $U=U(u_0)$ containing $\la$  and depending only on $u_0$, 
such that $$\dim \Ga_{u, U}=k.\tag 2.2$$
Furthermore,
$$\Ga_{u, U}=\Ga_{u_0, \la},\tag 2.3$$ as elements in $R(H)$.
\endproclaim

\demo{Proof} (2.2) is actually proved in~\cite{BW} (Lemma 17.1). 
More precisely, by~\cite{BW}, there exist a  $\e>0$ and $k$
continuous functions $$f_1,\ \dots ,\ f_k :(u_0-\e, u_0+\e)\to \R,\tag2.4$$ 
such that $$f_i(u_0)=\la;\tag 2.5$$ furthermore, for any $u\in (u_0-\e, u_0+\e)$, there exists an
open set $U$ that depends only on $u_0$ and contains $\la$ satisfying 
$$\{f_j(u)\}_{j=1}^k=\spec D_u \cap U.\tag 2.6$$

Let $Q_u$, $u_0-\e<u<u_0+\e$, be the orthonormal 
projections of $\Gamma_{u,U}$ onto $\Gamma_{u_0,U}$. By the 
continuity of $\{D_u\}$ and $f_j$'s, $\{Q_u\}$ is a 
continuous family of self-adjoint projections. Thus, it is possible to re-adjust $\e$ if necessary so that
$$\|Q_u-Q_{u_0}\|<1,\tag 2.7$$
for $u_0-\e<u<u_0+\e$. Now using a trick of Reed and Simon~\cite{RS, p.72},
if we define 
$$W_u=(1-(Q_u-Q_{u_0})^2)^{-{1\over 2}} [Q_u Q_{u_0}+ 
(1-Q_u)(1-Q_{u_0})],\tag 2.8$$
it is easy to verify that $W_u$ is unitary and 
$$W_u^{-1}Q_uW_u=Q_{u_0}.\tag 2.9$$
Notice that the image of $Q_u$ is $\Gamma_{u,U}$ and the above 
construction is $H-$compatible, (2.3) easily follows.
 \enddemo

Now we can proceed to define the equivariant spectral flow. Given  $\{D_u\}$ 
as above, we define $$\spec_H D_u=\{(\la, \Ga_{u,\la}); \la\in \spec 
D_u\}.\tag 2.10$$
By Lemma 2.2 and the fact that $R(H)$ has only countable many irreducible 
elements, there exist $f_j (u)\in C(I)$ and $R_j\in R(H)$, for $j\in \N$,
such that $$\cup_{u\in I} \spec D_u =\cup_j (f_j(u), R_j).\tag 2.11$$

Hence,  we can introduce the following:

\proclaim{Definition 2.3}
Given $D_u$ as above, we define the equivariant spectral flow as 
$$\sf_H(\{D_u\})=\sum_j \e(f_j) R_j,\tag 2.12$$
where $\e(f_j)$ is the intersection number of the graph $f_j$ with the 
line $u=-\delta$ for sufficiently small positive $\delta$.
We also denote
$$\sf(h, \{D_u\})=\tr(h, \sf_H(\{D_u\}).\tag 2.13$$
\endproclaim

\proclaim{Remark 2.4}
It is not hard to see that, as in the scalar case, only finite many
$\e(f_k)$'s in the above definition are non-zero. Also, if both $D_0$ and 
$D_1$ are invertible, $\delta$ can simply be replaced by $0$. \endproclaim

\proclaim {Remark 2.5} As in the scalar case, the equivariant spectral flow
is a homotopy invariant. In particular, let $E_\sigma$ be the affine space
of all the elliptic, positive-ordered operators in $\CH$ with the same
symbol $\sigma$. Then $E_\sigma$ is convex, and hence contractible.
Therefore given any two fixed $H$-compatible points in $E_\sigma$, the
equivariant spectral of two different $H$-compatible paths connecting them 
are the same. \endproclaim

\proclaim{Remark 2.6}Applying the method used above, it is not hard to 
extend the notion of higher spectral flow in the sense of Dai-Zhang~\cite{DZ}
to the equivariant setting. We leave the details to the 
interested readers.\endproclaim

For the rest of this paper, we pick a particular family $\{D_u\}$ as
below: 
$$D_u=(1-u)D+ug^{-1}D g,\tag 2.14$$ 
where $D$ is the Dirac operator given in (1.5) and being extended to the 
Hilbert space $\CH=L^2(S\otimes \C^N)$.

It is easy
to see that this family $\{D_u\}$ satisfy the conditions in Definition 2.1 
and 2.3.
Furthermore, we  note that $D_u$'s are of the same symbol, which is
in turn denoted as $\sigma$. 

Finally we want to prove  the following 
theorem to clarify the relation between
the equivariant spectral flow and our original index problem. 
The method used here is from Booss and Wojciechowski [BW]. 

\proclaim{Theorem 2.7}
$$\ind(h,T_g)=-\sf (h,\{D_u\}).\tag 2.15$$
\endproclaim

\demo{Proof}  Define $$P_u=(1-u)P+ug^{-1}P g,\tag
2.16$$ where $P$ is defined in (1.8). Apply the same argument used
in the proof of~\cite{BW} Theorem 17.17, noticing that it is compatible with
our equivariant setting. Thus, we conclude that
 $$\sf_H(\{D_u\})=\sf_H(\{P_u\}).\tag 2.17$$

Straightforward calculation gives that, 
 $$\ker(T_g)=\{u\in P_+ \CH, gu\in P_- \CH\},$$
$$\coker(T_g)=\{u\in P_- \CH, gu\in P_+ \CH\}.\tag 2.18$$

Then, it is easy to see that 

$$P_u(v)=\cases v, &  v\in P_+ \CH,\  gv\in P_+ \CH; \\
		 -v, & v\in P_- \CH,\  gv\in P_- \CH; \\
		(1-2u)v, & v\in\ker(T_g); \\
		(2u-1)v, & v\in\coker(T_g).
	\endcases  	\tag 2.19$$

Combining (2.17) and (2.19), (2.15) clearly follows.

\enddemo

\heading Section 3: Equivariant spectral flow 
and equivariant eta functions \endheading

Eta invariants first appeared  in [APS1], and has a known close 
relation with the 
spectral flow(cf.~\cite{BF, G}). In this section, we extend this 
relation to the equivariant case.
\proclaim{Definition 3.1}
Let $D$ be a self-adjoint operator on the Hilbert space $\CH$.
The eta function associated to  $D$ is defined to be
 
$$\eta (s,D)=\sum_{\la\neq0}{(sign\la){{\dim\Gamma_\la}\over {|\la 
|^s}}},\tag 3.1$$
where $Re(s)$  is large enough, $\la$ runs over the nonzero eigenvalues of 
$D$ 
and $\Gamma_\la$ is the eigenspace of $D$ with eigenvalue $\la$.
 
It is then clear that 
$$\eta (s, D)={1\over {\Gamma((s+1)/2)}} \int^{\infty}_0
{\tr(De^{-tD^2})t^{{(s-1)\over 2}}}dt \tag 3.2$$
holds. By a result of Bismut and Freed ([BF]), eta function of $D$ is 
analytic for $Re(s)>-1/2$, in particular, we write
 
$$\eta(D)=\eta(0,D).\tag 3.3$$
 
Furthermore, define the truncated $\eta$ function, for $\e>0$, to be
$$\eta_\e (s,D)={1\over {\Gamma((s+1)/2)}} \int^{\infty}_\e
{\tr(De^{-tD^2})t^{{(s-1)\over 2}}}dt\tag 3.4$$
and write
$$\eta_\e(D)=\eta_\e(0,D). \tag 3.5$$
\endproclaim

The equivariant eta function can be defined similarly:
\proclaim{Definition 3.2} Let $D$ be defined as in Definition 3.1. 
Furthermore, if there is compact group $H$ acting on $\CH$ and $D$  
commutes with the action of $H$,
the equivariant eta function associated to  $D$ is defined as
 
$$\eta (h,s, D)={1\over {\Gamma((s+1)/2)}} \int^{\infty}_0
{\tr(hDe^{-tD^2})t^{{(s-1)\over 2}}}dt, \tag 3.6$$
for $Re(s)$ large enough.
 
A regularity result of Zhang~\cite{Z} allows us to write
$$\eta(h,D)=\eta(h,0,D).\tag 3.7$$
 
We also define the truncated equivariant eta function, for an $\e>0$, to be
$$\eta_\e (h,s,D)=
{1\over {\Gamma((s+1)/2)}} \int^{\infty}_\e
{\tr(hDe^{-tD^2})t^{{(s-1)\over 2}}}dt,\tag 3.8$$
and
$$\eta_\e(h,D)=\eta_\e(h,0,D), \tag 3.9$$
for any $h\in H.$
\endproclaim

We then consider the variation of equivariant eta functions. 

Suppose $\F$ is the real Banach space of all bounded self-adjoint
operators on $\CH$. Let $\Phi$ be
the affine space $$\Phi=\{ D^S\otimes Id_{\C^N}+E, E\in \F\}.\tag 3.10$$ 
It is clear that for any $u$, $D_u$ as defined in (3.8) is in $\Phi$. 
\proclaim{Theorem 3.3}
For an $H$-invariant $D$ in $\Phi$ and any $h\in H$, define a one form 
$\alpha_{\e,h}$ on $\Phi$ such that for $X\in T_D \Phi=\F$,
$$\alpha_{\e,h}(X)(D)={(\e /\pi)^{\half} \tr(hX e^{-\e D^2})}. \tag 3.11$$ 
Then, $\alpha_{\e,h}$ is closed and we have
$$d\eta_\e(h,D)=2\alpha_{\e,h}(D).\tag 3.12$$
\endproclaim
\demo{proof}The proof of Theorem 3.4 is almost the same of that of Proposition 
2.5 of~\cite{G}, just noticing the commutativity of $h$ and $D$.
\enddemo

We can state  the main result of this section as:

\proclaim{Theorem 3.4} For  any $H$-invariant 
path in $\Phi=\Phi(D_0)$ connecting $D_0$ and $D_1$ as in (2.14), and any $h\in H$, we 
have 
$$\sf(h,\{D_u\})=-\int_\gamma \alpha_{\e,h} .\tag 3.13$$
\endproclaim
\demo{Proof} A similar formula for the scalar case is proved in [G]. Here 
we imitate the method used there.

By the fact that $\alpha_{\e, h}$ is closed and also Remark 2.5, both 
sides of (3.13)
are independent of the choice of the $H-$invariant path $\gamma: I\to 
E_\sigma$ with 
$\gamma(0)=D_0$ and $\gamma(1)=D_1$. $\cup \spec_H (\{\gamma(u)\}) $ 
can be written as $\cup_j(f_j(u),R_j)$ as in Section 2, where $R_j\in 
R(H)$ and $f_j\in C(I)$ for $j\in \roman{N}$.
 
Using a standard transversality argument, we can choose  an 
$H$-invariant path $\gamma$ such that for each $j$, the graph of $f_j$ 
intersects $\{u=0\}$ transversally. Also, by Remark 2.4, there is only 
finitely many nonzero $\e(f_j)$'s. Without loss of generality, let them 
be $f_1,\dots ,f_k$.

It is then easy to check that, for $h\in H$,
$$\sf(h,\{D_u\})=\sf(h,\gamma)=\sum_{j=1}^k \e 
(f_j)\tr(h,R_j).\tag 3.14$$

We can calculate the truncated equivariant eta function. 
For any $j\in\{1,\dots ,k\}$, the contribution of the $(f_j(u), R_j)$ to 
$\eta_\e(h,\gamma(u))$ for a $h\in H$, now denoted as  $S_{u,j}$, is
$${1\over \pi^{1/2}}\tr(h,  R_j)\int_\e ^\infty {f_j(u) e^{-t 
f_j(u)^2}t^{-1/2}dt}.\tag 3.15$$

Notice that 
$${1 \over \pi^{1/2}} \int_\e ^\infty {\la e^{-t
\la^ 2}t^{-1/2}dt} \to \pm 1, {\roman {if}}\ \la \to 0\pm.\tag 3.16$$
Hence, for $\tu$ being any zero of $f_j(u)$, let $\e(\tu)$ be the intersection 
number of $f_j(u)$ with $\{0\}\times I$ near $\tu$, we have
  $$S_{\tu+,j}-S_{\tu-,j}= 2\e(\tu)\tr(h,R_j).\tag 3.17$$
Summing up through all the zeros of $f_j$, and by the fact that
$$\sum_{\{\tu\in I;\ f_j(\tu)=0\}} \e(\tu)=\e(f_j),\tag 3.18$$
we have
$$\e (f_j)\tr(h, R_j)={1\over 2}\sum_{\tu;f_j(\tu)=0} 
(S_{\tu+,j}-S_{\tu-,j})$$
$$={1\over 2}((-\int_\gamma {dS_{u,j}}+S_{1,j}-S_{0,j}).\tag 3.19$$

Now summing up for all $j$, we are led to
$$\sum_{j=1}^k \e(f_j(u)) \tr(h,R_j)={1\over
2}(-\int_\gamma d\eta_\e(h,*)+\eta_\e(h,D_1)-\eta_\e(h,D_0)).\tag 3.20$$
  
Combine (3.12),(3.14) and (3.20), and notice that
$\eta_\e(h,D_1)=\eta_\e(h,D_0)$, we have (3.13).

\enddemo

Combining Theorem 2.7 and Theorem 3.4, we have the following:

\proclaim{Theorem 3.5} 

$$\ind(h,T_g)=\int_0^1{(\e /\pi)^{\half} \tr(h\dot{D_u}e^{-\e 
D_{u}^2})du}.\tag    3.21$$

\endproclaim

\proclaim{Remark 3.6}
Notice that the right-hand side of (3.21) is independent of the 
choice of $\e$, so we can use local index technique to calculate the limit 
of the integrand of (3.21) when $\e$ tends to 0. In such a way, we obtain a local version of Theorem 1.4.
\endproclaim

\heading Section 4: A Lefschetz theorem on odd spin manifolds  \endheading
In this section we apply the setting of~\cite{LYZ} to compute the right-hand 
side of (3.21).

First of all, we prove a Lichnerowitz type formula for $D_u^2$.

\proclaim {Lemma 4.1} We have
$$D_u^2=-\lap+{K\over 4}+u^2 c(\omega^2)+u(-\iota(\omega^*) \LC 
{S}+c(d\omega)+d^* \omega),\tag 4.1$$
where $$\omega=\dot D_u=g^{-1}dg,\tag 4.2$$
and  $K$ is the scalar curvature of $M$. \endproclaim
\demo{Proof} (4.1) follows easily from Prop 3.45 of~\cite{BGV} and the 
standard Lichnerowitz formula.
 \enddemo

For a fixed $h\in H$, let $F=\{x\in M;\ hx=x\}$  be the fixed point set 
of $h$. Without loss of generality we assume $F$ is a connected 
odd-dimensional totally geodesic sub-manifold and define its dimension to be 
$k$.  $\nu$ be the normal 
bundle of $F$ in $TM$, with dimension $2s$. Define $\nu(\delta)=\{x\in 
\nu;||x||<\delta\}$. Thus, $\nu$ to be invariant under $h_{TM}$, and $h_{TM}|\nu$ is 
non-degenerate.

If $P_\e(x,y):(S\otimes \C^n) \to (S\otimes \C^n)$ is the 
kernel for the operator $$O_\e =(\e /\pi)^{\half} 
\int_0^1 {(\dot{D_u}e^{-\e D_{u}^2})},\tag 4.3$$
 by the standard heat equation argument, we have 
$$\tr(h O_\e )=\int_M {\tr(h P_\e(hx,x)dvol}.\tag 4.4$$
A routine  argument  using  pseudo-differential operators shows that  if
$hx\neq x$, we have
$$\lim_{\e\to 0} \tr(h P_\e(hx,x))=  0.\tag 4.5$$

As a result, we may  localize the computation to $F$.

Notice that both
$\dot D_u$  and $h$ are bounded operators, working in a local 
trivialization of $\nu(\delta)$ yields the following:

\proclaim{Lemma 4.2}
Define, for $x\in F$,
$$L_{loc}(x)=\lim_{\e\to 0} \{\int_{\nu_\delta |x} {\tr(h P_\e(y, 
hy))dy}\}dvol_F,\tag 4.6$$ 
then it exists and is independent  of $\delta$.
Furthermore,
$$\ind(h, T_g)=\int_F{L_{loc}(x)}.\tag 4.7$$
\endproclaim

In the remaining part of this section, we  calculate $L_{loc}(x)$ for $x\in F$.

Fix any $x_0\in F$, let $e_1,\ \dots ,\ e_n\in TM $ be a local coordinate 
system in a neighborhood $\Cal N$ of $x_0$ such that $e_i$'s are 
orthonormal at $x_0$ and are parallel along
geodesics through $x_0$. And $e_1,\dots e_k\in TF$ and $e_{k+1},\dots ,
e_n\in \nu (F)$. For any $x\in \Cal N$ such that $hx\in \Cal N$, there is a  
$n\times n$-matrix $\FJ(x)$  satisfying

$$h|_{TM}(e_1(x),\dots ,\ e_n(x))=(e_1(hx),\dots ,\ e_n(hx))\FJ(x),\tag 4.8$$
while 
$$\FJ(x)=\pmatrix &1 & & & &\\
                & &\ddots& & &\\
                & & & 1 & &\\
                & & & & \exp(\Theta(x))\\
\endpmatrix, \tag 4.9$$
with $\Theta(x)\in so(2s)$.
 
Also,  denote $R^{TM}$ to be the curvature matrix of the Levi-Civita 
connection on $TM$ with respect to the chosen $\{e_i\}$:
$$(R^{TM})_{ij}=-{1\over 2}\sum_{p,q=1}^n {R_{ijpq}e^k e^l},\tag 4.10$$
for $1\leq i,\ j\leq n$, where $e^k$ is the dual vector of $e_k$. Also, 
if we choose the metrics and connections on $TF$ and $\nu(F)$ to be the 
restrictions of those of $TM$,respectively, we have the following curvature 
matrix

$$(R^{TF})_{ij}=-{1\over 2}\sum_{p,q=1}^k {R_{ijpq}e^p e^q}, \tag 4.11$$
for $1\leq i,\ j\leq k$ and
$$(R^{\nu(F)})_{ij}=-{1\over 2}\sum_{p,q=1}^k {R_{ijpq}e^p e^q},\tag 4.12$$
for $k+1\leq i,\ j\leq n$.

It is  known that $\FJ(x)$ is invariant along the fiber of $\nu$~\cite{BGV}. 
Hence, we can fix $e_i$'s such that we can have the following:

 $$\Theta(x_0)=\pmatrix &0&\theta_1&&&&\\
		         &-\th_1&0&&&&\\ 
			 &&& \ddots&&&\\
			 &&&&0&\th_s&\\
 			 &&&&-\th_s&0&\\
\endpmatrix,\tag 4.13$$	
with $0<\theta_i <2\pi$,

$$R^\nu(x_0)=\pmatrix  &0&v_1&&&&\\
                         &-v_1&0&&&&\\
                         &&& \ddots&&&\\
                         &&&&0&v_s &\\
                         &&&&-v_s  &0&\\
\endpmatrix,\tag 4.14$$

and

$$R^F(x_0)=\pmatrix      &0&&&&&&\\
			 &&0 &u_1&&&&\\
                         &&-u_1&0&&&&\\
                         &&&& \ddots&&&\\
                         &&&&&0&u_{(k-1)/2} &\\
                         &&&&&-u_{(k-1)/2}  &0&\\
\endpmatrix.\tag 4.15$$
Here all $u_i$'s and $v_i$'s are  two forms representing Chern roots.

It is easy to see that the kernel of $\sigma(h O_\e)$ is $\sigma(h P_\e)$, 
where the symbol map $\sigma$ is defined in (1.5).
  
Now we re-scale $T^*M$ as in [BGV] to get, for $\e\to 0$,
$$L_0=\lim_{\e\to 0} \sigma(h O_\e)$$
$$=\int_0^1{{h} ({1\over \pi})^{1/2} \omega \exp (\sum_i(\partial_i- 1/4 
\sum_j R^{TM}_{ij}b_j)^2 +u(1-u)\omega^2) du},\tag 4.16       $$
where $b_i$'s are local coordinate functions on $TM$ with respect to the 
chosen local charts.
 
We  proceed as in~\cite{LYZ} to get

$$Q_0(x_0,b)= \lim_{\e\to 0} {\sigma(h P_\e(x,hx))}$$
$$=\int_0^1  ({1\over \pi})^{1/2}[\omega \exp (u(1-u)\omega^2)] ({1\over 
{4\pi}})^{n/2}(-1)^s $$
$$\times (\prod _1^s \sin {\theta_i \over 2}) j_V(R^F) \exp({-1 \over 
4}\sum_1^s \sin(\theta_i) v_i (b_{k+2i-1}^2+b_{k+2i}^2))j_V(R^\nu) $$
$$ \times \exp(\sum_1^s(-\im {v_i\over 2}  \sin^2{\theta_i 
\over 2}\coth [\im {v_i\over 2}(b_{k+2i-1}^2+b_{k+2i}^2)]))du$$
 $$=\int_0^1 (-1)^s{1\over 2^n}  {1\over \pi^{m+1}}[\omega \exp 
(u(1-u)\omega^2)] j_V(R^F)(x)$$
$$\times \prod_1^s\{{(-1)\sin({\theta_i\over 2}) {{v_i/2}\over{\sin
v_i/2}}  \exp[- \sin ({\theta_i\over 2}){{v_i/2}\over{\sin v_i/2}}
\sin ({{v_i+\theta_i}\over 2}) (b_{k+2i-1}^2+b_{k+2i}^2)]}\}du
\tag 4.17$$
where
$$j_V(R^\nu)=\prod {{\im v_i/2}\over {\sinh (\im v_i/2)}},\tag 4.18$$
and
$$j_V(R^F)=\prod {{\im u_i/2}\over {\sinh (\im u_i/2)}}.\tag 4.19$$

In order to calculate $L_{loc}(x_0)$, we first apply (1.6) to  
$Q_0(x_0,b)$ and 
take the trace over $\C^N$;  then we integrate over $\nu(F)_{x_0}$. 
Asymptotically, that is to integrate over all 
$b_i$'s for $k+1\leq i \leq n$. Thus, we get $L_{loc}$ represented as a 
$k$-form on $F$ as following

$$L_{loc}(x_0)=({-\im\over{2\pi}})^{m+1}(\ch (g) j_V(R^F) \prod_1^s {[-\pi 
{\sin ({{v_i+\theta_i}\over 2})}^{-1}]})_{[k]},\tag 4.20$$
where $\ch (g) $ is defined as in (1.13).

Here notice that since every $\theta_i$ is non-zero, $\sin(\theta_i+v_i)^{-1}$ makes sense as a polynomial expansion.

It is then easy to get the characteristic class representation from 
(4.20), which is the following

\proclaim{Theorem 4.3} Let notations be as above. We have
$$L_{loc} (x_0)= ({-\im\over{2\pi}})^{m+1-s} (\hat A(F) \ch(g) 
[\pf (2\sin (\im(R^\nu +\Theta)/2)]^{-1} )_{[k]}(x_0). \tag 4.21
$$

\endproclaim
\proclaim{Remark 4.4} 
It is not hard to see that the method we have applied can also be used to 
prove similar local Lefschetz fixed point formulae for the Toeplitz 
operators associated to any Dirac-type operators.
\endproclaim

\enddocument 
\bye